%
%
%
%
%
\widowpenalty=15000
%
%
\outer\def\kopf #1 #2:#3\par{\medbreak
  \noindent{\bf(#1) #2:\enspace}{\it#3}\par
  \ifdim\lastskip<\medskipamount \removelastskip\penalty55\medskip\fi}
\abovedisplayskip=7pt plus 5pt minus 6pt
\belowdisplayskip=7pt plus 5pt minus 6pt
\hsize39pc \vsize52pc
\baselineskip11.5pt
\magnification=1200
\parindent 0pt

\def \ol{\overline}
\def \1{\backslash}

\def\exp{{\rm \hskip1.5pt  exp \hskip1.5pt }}


\font\ssf=cmss10

\font\calligX=callig15 at 10pt
\newfam\calligfam
\textfont\calligfam=\calligX
\font\calligVIII=callig15 at 8pt
\scriptfont\calligfam=\calligVIII
\font\calligV=callig15 at 5pt
\scriptscriptfont\calligfam=\calligV
\def\script{\fam\calligfam\calligX}
%
\font\eufmX=eufm10
\newfam\eufam
\textfont\eufam=\eufmX
\font\eufmVIII=eufm8
\scriptfont\eufam=\eufmVIII
\font\eufmV=eufm5
\scriptscriptfont\eufam=\eufmV
\def\frac{\fam\eufam\eufmX}
\font\eufbX=eufb10
\newfam\eubfam
\textfont\eubfam=\eufbX
\font\eufbVIII=eufb8
\scriptfont\eubfam=\eufbVIII
\font\eufbV=eufb5
\scriptscriptfont\eubfam=\eufbV

\def\ad #1{\hbox{$\raise1.75pt\hbox{$\scriptstyle /$}\!\hskip-0.5pt\hbox{\ssf A}_{ #1}$}}
\def\aq{\hbox{$\raise1.75pt\hbox{$\scriptstyle /$}\!\hskip-0.5pt\hbox{\ssf A}_{\bf Q}$}}
\def\parein#1#2{\par\noindent\rlap{\rm {#2}}\parindent#1\hang\indent\ignorespaces}
\def\paraus{\parindent 0pt\par}

\def\ar #1{\smash{\mathop{\pfr }\limits^{ #1}}}

\def\pfr{\hbox{$\hskip1pt\relbar\joinrel \hskip-1pt\relbar\joinrel
     \hskip-2pt\longrightarrow\hskip1pt$}}

\def\gmm #1{\hbox{$({\bf G}_m)_{{#1}}$}}
\def\gma #1 #2{\hbox{$({\bf G}_m)_{{#1}}({#2})$}}
\def\gm1 #1{\hbox{$({\bf G}_m^1)_{{#1}}$}}
\def\h1 #1 #2{\hbox{$H^1({#1}, {#2})$}}
\def\hn #1 #2{\hbox{$H^n({#1}, {#2})$}}
\def\rs #1 #2 #3{\hbox{${\rm R}_{{#1}/{#2}}({#3})$}}
\def\rg #1 #2 #3{\hbox{${\rm R}_{{#1}/{#2}}\,\gmm {#3}$}}
\def\r1 #1 #2 #3{\hbox{${\rm R}_{{#1}/{#2}}(\gm1 {#3})$}}

\def\no #1 #2 #3{\hbox{$ N_{{#1}/{#2}}({#3})$}}
\def\noabb #1 #2 {\hbox{$ N_{{#1}/ {#2}}$}}
\def\mmatrix#1{\null\,\vcenter{\normalbaselines\mathsurround=0pt
    \ialign{\hfil$##$\hfil&&\hskip1pt\hfil$##$\hfil\crcr
      \mathstrut\crcr\noalign{\kern-\baselineskip}
      #1\crcr\mathstrut\crcr\noalign{\kern-\baselineskip}}}\,}
\def\epmatrix#1{\left(\null\vcenter{\normalbaselines\mathsurround=0pt
    \ialign{\hfil$##$\hfil&&\enspace\hfil$##$\hfil\crcr
      \mathstrut\crcr\noalign{\kern-\baselineskip}
      #1\crcr\mathstrut\crcr\noalign{\kern-\baselineskip}}}\right)}
\def\tpmatrix#1{\left(\null\vcenter{\normalbaselines\mathsurround=0pt
    \ialign{\hfil$##$\hfil&&\thinspace\hfil$##$\hfil\crcr
      \mathstrut\crcr\noalign{\kern-\baselineskip}
      #1\crcr\mathstrut\crcr\noalign{\kern-\baselineskip}}}\right)}
\def\ematrix#1{\null\vcenter{\normalbaselines\mathsurround=0pt
    \ialign{\hfil$##$\hfil&&\enspace\hfil$##$\hfil\crcr
      \mathstrut\crcr\noalign{\kern-\baselineskip}
      #1\crcr\mathstrut\crcr\noalign{\kern-\baselineskip}}}}
\def\tmatrix#1{\null\vcenter{\normalbaselines\mathsurround=0pt
    \ialign{\hfil$##$\hfil&&\thinspace\hfil$##$\hfil\crcr
      \mathstrut\crcr\noalign{\kern-\baselineskip}
      #1\crcr\mathstrut\crcr\noalign{\kern-\baselineskip}}}}
\def\gal #1 #2{\hbox{${\rm Gal}({#1}\vert {#2})$}}
\def\agal #1 {\hbox{${\rm Gal}(\ol{#1}\vert {#1})$}}


\def\xtyz #1 #2 #3{\hbox{$ {#1}\otimes_{{#2}}  {#3}$}}

\def\deff{{\buildrel {\rm def} \over =} }

\def\fa{\c F_{\rm aut}}

\def\>{\rangle}
\def\<{\langle}

\abovedisplayskip=8pt plus 5pt minus 3pt
\belowdisplayskip=8pt plus 5pt minus 3pt

\font\ninecaps=cmcsc9
\font\ninebf=cmbx9
\font\nineit=cmti9
\font\ninerm=cmr9

\font\eightbf=cmbx8
 at 11truept

\def\re{{\rm Re}\,}

\def\la{{\script L}\hskip4.7pt}
\def\fh{{\script F}\hskip4.7pt}

\def\hhfactor{\nu}
\def\HHfactor{N}

\def\erfc{{\rm Erfc}\,}

\def\law{{\buildrel {\rm law} \over =} }
%
%
%
%
%
%
\output={
\ifodd\pageno\hoffset=0pt
\else
\hoffset=0truecm\fi\plainoutput}
\font\authorfont=cmcsc10 at 11pt
\font\titlefont=ptmb at 11pt

\font\headingfont=cmbx12 at 10pt
\font\paragraphfont=cmbx10

\font\ninecaps=cmcsc9
\font\ninebf=cmbx9
\font\nineit=cmti9
\font\ninerm=cmr9

\font\eightbf=cmbx8

\def\re{{\rm Re}\,}

\def\la{{\script L}\hskip4.7pt}
\def\fa{{\script F}\hskip4.7pt}
\def\ndlint #1 #2{\hbox{$\, \int\limits_{#1}^{#2}\,$}}
\def\ndint #1 #2{\hbox{$\,\int_{#1}^{#2}\,$}}
\def\ppar{\vskip2.5pt plus1pt minus1pt}

\abovedisplayskip=6pt plus 5pt minus 3pt
\belowdisplayskip=6pt plus 5pt minus 3pt
\bigskipamount=11pt plus 3pt minus 3pt       
\medskipamount=5pt plus 1.5pt minus 1.5pt
\smallskipamount=2.5pt plus 1pt minus 1pt
\baselineskip=12pt
%
\voffset.45cm
\centerline{\titlefont ON A BROWNIAN EXCURSION LAW, I:}
\centerline{\titlefont CONVOLUTION REPRESENTATIONS}
\baselineskip=11.5pt
\vskip.3cm
\centerline{\authorfont Michael  Schr\"oder}
\vskip.3cm
\centerline{
\vbox{
\hsize=10.8cm
\baselineskip=9.3pt
\ninerm
{\ninebf Abstract:} This paper studies Brownian motion subject to the
occurrence of a minimal length excursion below a given excursion level.
The law of this process is determined.  The characterization is explicit 
and shows by a layer construction how the law is built up over time in 
terms of the laws of sums of a given set of independent random variables.
\vskip5pt
{\ninebf Key words:} Brownian excursions, Brownian law subject to
excursion conditions,  Laplace transform and its inversion. 
\vskip5pt
{\ninebf MSC 2010:} Primary  60J65, 60G46.}}
\vskip9pt
\centerline{\headingfont 1.\quad Introduction}
\vskip.1cm
In this paper we identify a general Parisian-type excursion law and 
determine its structure in the Brownian case. The idea is that we do 
not just ask when a process passes a given level, but instead ask 
for how long it will thereafter stay on the same side of the level. 
These excursions may have any length and may occur at any time. Having 
settled on a minimal length for them, the question to ask is if a 
minimal-length excursion will occur during a given time span. This 
then yields a second source of randomness for the problem, and the 
object to be studied is its joint law with that of the given process; 
see Section~2 for how to make all this precise. 
\ppar
The principal difficulty is to obtain an explicit description 
of this joint law when the process to start with is given  explicitly. 
The main contribution of the paper is to give such an explicit
description when this process is Brownian motion; details are in
Section~3. 
\ppar
Historically, a rigorous study of these questions seems
to originate with [13]. Extensions and further
applications are developed for example in [6] and the
references there. Applications to insurance are emerging; see for 
example [5]. 
\ppar
Regarding the structural understanding
of the joint law, its Laplace transform with respect to time has 
been determined in the Brownian case; see
[13, Appendix] with modifications in [11].
As made precise in Section~5, the Laplace transform permits one to
separate the above two sources of randomness
of the problem. This is effected in [13, Appendix] by recourse
to the Brownian meander and the Az\'ema
martingale; the key steps of the argument are recalled in Section~6. The
result is quotients of higher transcendental functions
as Laplace transform of the joint law, with the denominators
corresponding to the law's excursion source of stochasticity  and the
numerators to the law's process source of stochasticity. This generalizes
situations as encountered for example with Laplace
transforms of first passage times of Brownian motion, which lead to
laws in terms of theta functions.
\ppar
Mathematically, the contribution of the paper is to provide analogous
such representations to the present context. The basis for this is furnished 
by our finding in Appendix~B that the above denominators satisfy a functional
equation, and as its main technical contribution the paper turns this
into a way for the analytic inversion of the corresponding quotients.
With more details developed in Section~4, this then shows how the effect 
of the excursion source on the law is built up over time --  by adding at any
integer time $n\!+\!1$ a new layer of $n$ independent copies of a
fixed explicit random variable. This picture is made rigorous
in Section~9 yielding the explicit expressions for the joint law of 
Section~3.
\ppar
To conclude,  Section~10 indicates how our approach is 
complemented in the companion [12] to the present paper 
by constructive techniques for convolution representations.   
\vskip8pt
\centerline{\headingfont 2.\quad Basic notions and facts:
the general process excursion law}
\vskip.08cm
We address in the paper the occurrence of
excursions of a minimal length duration of a given stochastic  process below
a given excursion level which take place within a given period of time,
and ask to characterize the law of the new stochastic process thus
obtained. We follow [11] to make this precise in two stages. 
A general framework is thus established in the present section which 
furnishes the basis for the Brownian case to be developed subsequently 
in the paper.
\medskip
{\bf 2.1\quad Basic setting:}\quad Let $X$ be any real-valued
stochastic process on the complete probability space $(\Omega, {\cal F}, Q)$
whose time set is the nonnegative reals $[0,\infty)$, so that
$X=(X_u)_{u\ge 0}$, and let it be adapted to the filtration 
${\bf F}=({\fh}_u)_{u\ge0}$
on this space which satisfies the usual conditions. 
\ppar
As the primitives of the problem fix any time $t$,
specify any real $a$, which is to take the r\^ole of the {\it excursion
level\/} of the problem, and specify any real $D>0$, which is  to take
the r\^ole of the {\it minimal excursion duration\/}.
\ppar
Denote by $T_{t,a}:=T(X)_{t,a}$ the first passage time after time $t$ 
to the level $a$ of $X$ given by $T(X)_{t,a}=\inf\{ s\ge t| X_s\in [a,\infty)\}$; 
for simplicity assume that $T_{t,a}$ is a stopping time, and 
if $t=0$ we drop reference to the time $t$ from the notation. 
\medskip
{\bf 2.2\quad Achievement time:}\quad We first recall how the excursion 
situation under consideration is  encoded by the {\it achievement time\/}
$H_{a,t}:=H(X)_{a,t}$, a random variable to be defined presently with 
the property that for any real $T>t$ we have $H_{t,a}\le T$ iff there is in
$[t,T]$ a subinterval of length at least $D$ on which $X$ only takes values
below the level $a$. Morally, $H_{a,t}$ so is the upper bound of 
the first of these subintervals occurring after time $t$. For its 
explicit construction, which is to follow, for any real $u\ge 0$ 
introduce $\sigma_{a,u}=\sigma(X)_{a,u}$ to denote the smallest 
upper bound of subintervals as above of $[u,\infty)$, namely:
$$
\sigma_{a,u}=\inf\{ s\in[u\!+\!D,\infty)\,|\, X|_{(s-D,s)}<a\}\,.\quad
u\in[0,\infty)\,.
$$
We will need to distinguish 
the two cases when
$X_t$ is equal to or above the level $a$, and 
when it is below this level.
\goodbreak\ppar
If $X_t\ge a$,  the situation to be referred to as the {\it first case\/}, 
the question arises if $X$ will become and remain smaller than $a$ for all 
points in time of an interval of length at least $D$ contained in $[t,\infty]$. 
We represent this first case by setting:
$$H_{a,t}=\sigma_{a,t}\quad \hbox{if $X_t\ge a$}\,.$$
If $X_t<a$, the situation to be referred to as the {\it second case\/}, 
a new idea is needed when $X$ has been staying below $a$ for some connected 
period of time already. This excursion needs to continue below
$a$ for a period  of length  $\delta_t \in (0,D)$ only to reach the 
minimum duration $D$ stipulated.  Two subcases result which we formalize 
in terms of the location of $t\!+\!\delta_t$ relative to the first passage 
time $T_{t,a}=T(X)_{t,a}$. 
Firstly, if $T_{t,a}> t\!+\!\delta_t$, the situation to be referred to as 
{\it Subcase~1\/}, then the process $X$ continues to stay below $a$ during 
the period of time from time $t$ until time $t\!+\!\delta_t$ a fortiori,
whence the representation:
$$H_{a,t}=t\!+\!\delta_t\quad \hbox{if $X_t<a$ and $T_{t,a}> t\!+\!\delta_t$}\,.$$
Secondly, if $T_{t,a}\le t\!+\!\delta_t$, the situation to be referred to as 
{\it Subcase~2\/}, the level $a$ is passed by $X$ earlier than time $t\!+\!\delta_t$,
and the clock for reaching the minimum length $D$ is restarted at time
$T_{t,a}$. This then puts us into the first case, albeit as of time $T_{t,a}$
instead of time $t$, whence the representation:
$$H_{a,t}=\sigma_{a,T_{t,a}}\quad \hbox{if $X_t<a$ and $T_{t,a}\le t\!+\! \delta_t$}\,.$$
{\bf 2.3\quad Excursion density:}\quad 
The excursion law to be studied is then by definition the representing 
measure of the functional given by:
$$
\phi\mapsto E\big[{\bf 1}_{\{H_{a,t}\le T\}}\, \phi(X_T)\big|\, \fh_t\big]\,,
$$
for any  $Q$-measurable function $\phi$ on the reals {\bf R}; assuming 
absolute continuity with respect to Lebesgue measure in addition the 
concept to be studied is hence the {\it  excursion density\/} 
$h_{a,t}:=h(X)_{a,t}$ which by definition is the conditional density on  
$[t,\infty)\times {\bf R}$~characterized~by:
$$
E\big[{\bf 1}_{\{H_{a,t}\le T\}}\, \phi(X_T)\big|\, \fh_t\big]
=\ndint {\bf R} { }    \phi(x)\, h_{a,t}(T,x)\, dx\,,
$$
for any  $Q$-measurable function $\phi$ on the reals {\bf R}; 
here notice $x\mapsto h_{a,t}(T,x)=0$ for $T<\delta_t$ 
extending the definition of $\delta_t$ in Section~2.2 to the first 
case there by the convention $\delta_t:=D$.
\medskip
{\bf 2.4\quad Normalizations:}\quad We develop a normalized form of the excursion 
law. The normalization is based on an emulation of the re\-starting at a fixed 
stopping of Markov processes. The construction starts from postulating the 
existence of a process $X^*$ which satisfies 
$$
X^*_u\law X_{t+u}\!-\!X_t\, , \qquad u\in[0,\infty)
\,  ,$$
as well as $X^*(0)=0$, and which is independent of $\fh _t$; examples 
for processes $X$ for which these normalized processes exist thus include
Brownian motion, L\'evy processes which creep across levels, or 
continuous Markov processes which are homogeneous in their time and their space 
variable. The construction of Section~2.2 applied to $X^*$ then yields 
for any real $b$  the 
{\it normalized achievement time\/} $H_b^*:=H(X^*)_{b,0}$ which satisfies 
$$\displaylines{
H_{a,t}=t\!+\! H_{a-X_t}^*\,,\cr  
\noalign{as well as}
E[{\bf 1}_{\{H_{a,t}\le u+t\}}\, \phi(X_{u+t})|\fh_t]
=E[{\bf 1}_{\{H^*_{a-X_t}\le u\}}\, \phi(X_t\!+\!X^*_u)]\,,
\cr}$$
for any measurable map $\phi$ on {\bf R}, 
treating the time-$t$ value $X_t$ of $X$ as a real number.
Granting absolute continuity with respect to Lebesgue measure in addition 
as in Section~2.3 above, densities $h_b^*$ on $[0,\infty)\times {\bf R}$ 
are hence defined by:
$$
E\big[{\bf 1}_{\{H^*_b\le u\}}\, \phi^*(X^*_u)\big]
=\ndint {\bf R} {} \phi^*(y) h^*_b(u,y)\, dy\,,
$$
for any $Q$-measurable function $\phi^*$ on  {\bf R}. By construction
these {\it normalized excursion densities\/} $h^*_b$ satisfy 
$h_b^*( [0,\delta_0)\times {\bf R})=\{0\}$ with
$\delta_0$ the minimal excursion duration remaining
as in Section~2.3 above, and 
are related to the excursion densities $h_{a,t}$ of Section~2.3 by:
$$h_{a,t}(u\!+\!t,y\!+\!X_t) =h^*_{a-X_t}(u,y)\,,$$
for any real $u\ge 0$ and $y$. Granting existence, the study of the 
former law is hence reduced to the latter one. Translating 
the first and second case of Section~2.2 this explicitly
asks to characterize $h_b^*$, for any real $b$, in
the following
{\it Cases I and II\/} respectively:
\medskip
\parein{25pt}{(I)} We have $b\le 0$ and ask if there will be a length
$D$ time interval $I_D$ on which $X^*|_{I_D}<b$.
\paraus
\parein{25pt}{(II)} We have  $b>0$ and ask if
$X^*|_{[0,\delta_0)}<b$ where the positive real $\delta_0<D$ is the
Section~2.2 minimal excursion length remaining.
\paraus
\medskip
Here Case~II further decomposes into two subcases according
to the relative position of  $\delta_0$ to $T_b^*$, the first passage
time of $X^*$ to the level $b$, as follows. If moreover  
$T_b^*>\delta_0$, we are in subcase (II-1), otherwise we are in 
subcase (II-2). 
\ppar
The paper initiates a study of the structure of $h_b^*$ along these lines 
when $X$ is Brownian motion; our main results here are to 
be described in the next Section~3 as a first step.
\vskip8pt
\centerline{\headingfont 3.\quad Statement and discussion of main results:}
\centerline{\headingfont \phantom{4.\quad}  convolution
representations of the excursion law in the Brownian case}
\vskip.08cm
In this section we formulate and discuss our main results about the
normalized excursion law $h_b^*$ of Section~2.4 in the Brownian
case. We delineate in Sections~3.2 and 3.3 how, with hindsight, it 
emerges from a core relationship pinned down in Section~4 by 
taking into account the additional interrelations entailed by the 
Case~I and Case~II situations. Our discussion is in terms of 
convolutions of laws, 
with representing functions listed in Section~3.1 to follow.  
\medskip
{\bf 3.1\enspace Set up and functions:}\quad 
We modify the Section~2.1 set up as follows. 
First we let the probability space $(\Omega, {\cal F}, Q)$ 
there be equipped with the  standard filtration 
${\bf F}=(\fa _u\,|\,u\in [0,\infty))$ of the
Brownian motion $W$ on it. Then let $X=W$ and, consequently, let $X^*=W^*$ 
in Section~2.4 with $W^*$ now being the standard Brownian motion independent 
of time-$t$ information $\fa_t$ that is obtained by time-$t$ restarting of $W$.
\ppar
It is convolution representations for $h_b^*$ which thus result. They 
are in terms of three main classes of functions. Referring to Appendix~A 
for more detail about them, these are as follows. 
Firstly, the functions $\hhfactor_n$ on $(0,\infty)$  given for any integer 
$n\ge 1$ by the $n$-fold convolutions on $(0,\infty)$: 
$$ \hhfactor_n=\hhfactor^{*(n)}\quad \hbox{where}\quad 
\hhfactor(u)=(2/\sqrt{\pi})\sqrt{u}/(2u\!+\!1),
\quad u\in [0,\infty)\,,$$
and by $\hhfactor_0=\hhfactor^{(0)}$ being set equal to the Dirac delta 
function at $0$; for $n\ge 1$ note $\hhfactor_n={\rm law}\, (nV)$ with $nV$ the 
$n$-fold sum of independent copies of a random variable $V$ with 
${\rm law}( V)=\hhfactor$.  Secondly, for any real $a$ and $c$, the
functions $\rho_{a,c}$ on $(0,\infty)$ given by:
$$\rho_{a,c}=\cases{\sqrt{D}\, \chi_{|a+c|}*\hhfactor& \quad $c\le 0,$\cr
\noalign{\vskip1.5pt}
\sqrt{D}\,  \chi_{|a|}*g_c& \quad
$c>0;$\cr}
$$
they are convolutions on $(0,\infty)$ in terms of 
two additional functions on $(0,\infty)$ that depend  on the 
complex parameter $\alpha$ with $|\arg(\alpha)|\le {\pi\over 4}$ 
and for any real $u\ge 0$ are given  by: 
$\chi_\alpha(u)=(1/\sqrt{\pi u}\,)\exp(-(\alpha/2)^2/u)$
respectively:
$$
g_\alpha(u)=\exp\Big( \!-{\alpha^2/2\over 1\!+\!2u}\Big)
\Big\{ h(u)\exp\Big( \!-{\alpha^2\over 4u(1\!+\!2u)}\Big) +
{\alpha\over 1\!+\!2u}
\erfc\Big( {\alpha \over \sqrt{4u(1\!+\!2u)} }\Big)\Big\}.
$$
Thirdly, for any  excursion level $b$ the function $\beta=\beta_b$
on the real line {\bf R} given~by:
$$
\beta(y)=(b\!-\!y)/\sqrt{D}\,.$$
{\bf 3.2\quad Case I results:}\quad We first assume the situation of 
Case~I to hold as it is described in Section~2.4; we are thus asking 
about the occurrence of an excursion of duration at least $D$ below 
the level $b\le 0$ at some point in the future. The structure of the 
normalized excursion law $h^*_b$ here is governed by the sign of 
$\beta$ on the `state space' by way of the functions $\rho_{a,c}$. 
Our precise result is as  follows.
\medskip
{\bf Theorem 3.1:}\quad {\it In the above Case~{\rm I} setting
we have for any reals $u>0$ and $y$  the convolution
sum representation:
$$
h^*_b(u,y)={1\over 2D}\sum\nolimits_{1\le n<{u\over D}}
{(-1)^{n-1}\over (2\pi)^{n/2}}\,
\Big(\rho_{b^* , \beta(y)}*\hhfactor^{*(n-1)}\Big)
\Big( (1/ 2)\big({u\over D}-n\big) \Big),
$$
where the sum is over all integers $n\ge 1$ satisfying $n<u/D$,
and where $b^*=b/\sqrt{D}$.}
\goodbreak\medskip
With a proof in Section~9.1 
our results make explicit how the law is
built up over time by the addition of a new layer at integer points in
time. Anticipating the discussion in Sections~5 and 6 this is a principal
feature. It gives expression to the effect of $X_1=H_b^*$ as the first of 
two sources of stochasticity that determine the structure of the excursion 
law. The effect 
of the second source of stochasticity $X_2=W^*(H_b^*)$, on the other
hand, is encoded by the structure of the layers to be added. Here it gives
rise to a fixed density $\rho=\rho_{b^*,\beta(y)}$ which, at each time $n$, 
is convolved with $\hhfactor_{n-1}$, the $(n\!-\!1)$-fold convolution of 
the density $\hhfactor$ originating with the above first source of 
stochasticity. The structure of the function $\rho$ depends on the state 
variable $y$, and gives expression to the  position of $y$ relative to the 
excursion level $b$ by way of the function~$\beta$.
\medskip
{\bf 3.3\quad Case II results:}\quad We next 
address the structure of the normalized excursion law $h_b^*$ in the
situation of Case~II of Section~2.4; recall that this setting 
is characteristic for the problem by asking about the excursion of Brownian 
motion $W^*$ to continue below the level $b>0$ for a period of length at 
least $\delta:=\delta_0>0$ smaller than $D$.  We here obtain
$h_b^*$  as  a four~term~sum:
$$\hbox{$
h^*_b=\sum\nolimits_{j=1}^2 \sum\nolimits_{k=1}^2 h^*_{b,j,k}
$}\,,
$$
where the single summands are functions of the time variable $u>0$
and the state variable $y$ as follows. Pertaining to the
situations where the  first source of stochasticity $X_1=H_b^*$ 
(see the discussion of Section 3.2) is in Case~II-1, 
$$\displaylines{
h^*_{b,1,1}(u,y)={Q^*(T_b>\delta)Q^*(T_b\le \delta)\over D}\,
{{\bf 1}_{(\delta,\infty)}(u)\over \sqrt{2\pi(u\!-\!\delta)}}
\ndint 0 {\infty}  x\exp\Big(\!-{x^2\over 2D}
-{(b\!-\!x\!-\!y)^2\over 2(u\!-\!\delta)}\Big)dx,\cr
\noalign{\vskip1.5pt}
h^*_{b,1,2}(u,y)={Q^*(T_b>\delta)}\,
{\bf 1}_{(\delta,\infty)}(u)\varphi_{b,u}(y);
\cr}$$
here $\varphi_{b,u}$ is a translate of the transition
density from $0$ to $y<b$ of Brownian motion killed at the first hitting
time of $b$ and living on $(-\infty, b)$ as follows:
$$\varphi_{b,u}(y)={1\over \sqrt{2\pi u}}\exp\Big(\! -{y^2 \over 2u}\Big)
- {1\over \sqrt{2\pi u}}\exp\Big(\! -{(y\!-\!2b)^2 \over 2u}\Big)\,,$$
where $Q^*(T_b>\delta)=(1\!-\!\erfc)(b/\sqrt{2\delta}\,)$ and
$ Q^*(T_b\le \delta)=\erfc (b/\sqrt{2\delta}\,)$, and where the integral
factor of $h_{b,1,1}^*$ is expressible in terms of $\erfc$ and its derivative
as well; see the representation for the function $h_{b,3}$ in [11, p.4].
\ppar
The remaining two summands have been known  at the level 
of their Laplace
transforms. By inverting these transforms in Sections~9.2 and 9.3 we are 
now able to describe in the next two results how their structure is built up
over time. This once more proceeds by expressing them as sums of convolutions
in terms of Section~3.1 functions. Our first result pertains to the
summand $h^*_{b,2,1}$ as follows.
\goodbreak\medskip
{\bf Theorem 3.2:}\quad {\it In the above Case~{\rm II} situation, with
excursion level $b>0$, we have for any reals $u>0$ and $y$ the
convolution sum representation:
$$
h^*_{b,2,1}(u,y)={1\over 2D}\sum_n
{(-1)^{n-1}\over (2\pi)^{n/2}}\ndint 0  {\delta\wedge (u-Dn)}
\hskip-3pt\mu_b(dw)\,
\Big(\rho_b(\cdot,y)*\hhfactor^{*(n-1)}\Big)\Big( {u\!-\!Dn\!-\!w\over 2D}\Big),
$$
where the sum is over all integers $n\ge 1 $ satisfying $n<u/D$.}
\medskip
With a proof in Section~9.2 this theorem refers to a situation where the
first stochasticity source $X_1=H_b^*$ is now considered on the set of events 
where is expected to be restarted. This is given expression to by integration 
with respect to $\mu_b$, the law  of the 
first passage time of Brownian motion to the level $b$ given by: 
$$
\mu_b(dw)=\psi_{b\sqrt{2}}(w)\, dw\,;
$$
see Section~A.1. Conditional on that influence, the excursion law  is
built up at positive integer time points as follows. At each time $n$
add the new term obtained by convolving  $\hhfactor_{n-1}$, the 
$(n\!-\!1)$-fold convolution of the density $\hhfactor$, with a fixed
density $\rho_b(\cdot,y)$. Here the function $\rho_b$ on
$(0,\infty)\times {\bf R}$ is explicitly given by the
following integral:
$$
\rho_b(\tau,y)=2\ndint {\bf R} {}  \varphi_{b,\delta}(x)\, N_{0,D\tau} \big(
|x\!-\!y|/\sqrt{2}\,\big)\, dx\,,
$$
where $N_{0,v}(\xi)
=(1/\sqrt{2\pi v})\int_{\raise2pt\hbox{$\scriptstyle (-\infty,\xi]$}}
\exp(-x^2/(2v))\, dx$
denotes the normal distribution with mean equal to $0$ and variance
equal to $v$. The functions $\rho_b$ originate with the restriction of the 
second source of stochasticity $X_2=W^*(H_b^*)$ to the set of all events 
where a restarting of $H_b^*$ will not happen. 
\medskip
The situation for the final fourth summand  $h_{b,2,2}^*$ differs from
the above set up in regard to the stochasticity source $X_2=W^*(H_b^*)$: 
here it is considered on the set of all events where a restarting  of its
argument $H_b^*$ is to take place. Conditional on the restarting of 
$X_1$ as above we therefore are in the situation of Section~3.2, and 
hence seek to characterize the occurrence of future length $D$ excursions 
below the level $0$. This situation is encoded by the Section~3.1 functions
$\rho_{0,\beta(y)}$ and  our precise result is as follows.
\medskip
{\bf Theorem 3.3:}\quad {\it In the above Case~{\rm II} situation, with
excursion level $b>0$, we have for any reals $u>0$ and $y$ the
convolution sum representation:
$$\eqalign{
&h^*_{b,2,2}(u,y)\cr &={Q^*(T_b\le \delta)\over 2D}\sum_n
{(-1)^{n-1}\over
(2\pi)^{n/2}}\!\ndint 0  {\delta\wedge(u-Dn)} \hskip-4pt\mu_b(dw)\,
\Big(\rho_{0,\beta(y)}*\hhfactor^{*(n-1)}\Big)\Big( {u\!-\!Dn\!-w\over 2D}\Big),\cr}
$$
where the sum is over all integers $n\ge1$ satisfying $n<u/D$.}
\medskip
%
This result is proved  Section~9.3,  as the final
step of our argument which starts in Section~4.
\vskip5pt
\centerline{\headingfont 4.\quad The key result for Laplace inversion}
\vskip.08cm
At the heart of the structure  of the normalized excursion law 
$h_b^*$ as expressed by the results of Section~3 is the reconstruction
of the law of $H_0^*$, the level-$0$ normalized achievement time. This
reconstruction is to be established in this section, and asserts 
$H_0^*$ morally to originate as (the weak limit of) an infinite sum of 
independent copies of a fixed random variable, to be denoted by $V$, 
in which only finitely many summands contribute when looked at in a 
pointwise sense. Here $V$ is pinned down in terms of its distribution 
by:
$${\rm law} ( V) =\hhfactor
\quad\hbox{where}\quad    
\hhfactor(u)=(2/\sqrt{\pi}\,)\, \sqrt{u}/(2u\!+\!1)\,,\quad u\in [0,\infty);
$$
our construction, however, necessitates two auxiliary
random variable summands as follows. 
\medskip
{\bf Theorem 4.1:}\quad {\it For any real $u>0$, we have the equality of 
laws of pairwise independent random variables:
$$
{\rm law} \big( A_R\!+\!H_0^*\big)(u)  
=
\sum\nolimits_{1\le n<2n}\, {(-1)^{n-1}\over (2\pi)^{n/2}}\, 
{\rm law} \big( B_R\!+\! (n\!-\!1)V\big) \big( u\!-\hbox{${1\over 2}$}n\big)\,,
$$
where $A_R$ and $B_R$ are random variables such that
on some  complex half-plane $\{ \re z >z_0\}$ with $z_0\ge 0$ we have:
$$E[\exp(-zA_R)]=R(z) 
\quad \hbox{and}\quad 
E[\exp(-zB_R)]=R(z)/\sqrt{z}\,,
$$
for an analytic function $R$ which for some $a>{1\over 2}$
satisfies $R(z)=O(|z|^{-a})$ as $|z|\to\infty$.}
\medskip
The starting point of this result is provided 
by joining the equality
$$E[\exp(-zH_0^*)]=1/\Psi(\sqrt{2Dz}\,)\,, \quad \re(z)>0\,,
$$
anticipated from Section~6.4  in terms of the  function 
$\Psi(w)=\int_{\raise 2pt\hbox{$\scriptstyle (0,\infty)$}}
x\exp(-x^2/2\!+\!wx)\, dx$ of Appendix~B,  with the symmetry 
property of $H_0^*$ expressed by the functional equation
$$
\Psi(w)=\Psi(-w)+\sqrt{2\pi}\,w\exp(\hbox{${1\over 2}$}w^2)\,,
$$
for any complex $w$, furnished by the key identity of Appendix~B.
This enables a representation of the expectation by resolution as 
a geometric series:
$$ {R(z)\over \Psi(\sqrt{z\, }\, )}
=R(z){f(z)\over 1\!-\!p(z)}
={R(z)\over \sqrt{2\pi z\, }}\sum\nolimits_{n=0}^\infty {(-1)^{\scriptstyle n}
\over (2\pi)^{\scriptstyle n/2}}
\, \exp\big(-\hbox{${1\over 2}$}(n\!+\!1) z\big)N_1(z)^n,$$
where $f(z)=(1/\sqrt{2\pi z}\,) e^{-z/2}$ 
and $p(z)=-f(z)\Psi(-\sqrt{z}\, )$, where $N_1(z)=\Psi(-\sqrt{z}\,)/\sqrt{z}$
from Appendix~A, and granting that for some $a_0>0$  the first two functions 
are smaller than $1$ in absolute value when $\re(z)\ge a_0$. Granting too that 
Laplace inversion of this sum can then be effected term by term, 
the assertion of Theorem~4.1 follows noting $\la^{-1}(N_1^n)=\hhfactor_n$ 
from Proposition~A.1. With the analytical facts implied by 
the Appendix~B leading term expansion of $\Psi$, omitting further detail 
the proof of Theorem~4.1  is complete.  
\medskip
{\bf Remark 4.2:}\enspace The proof of Theorem~4.1 can be seen
as being effected by transfer to the framework of the It\^o theory of excursions
of Brownian motion. Denoting the objects of this theory  by $e$, this
perspective  then in particular affords an interpretation of the 
function $p$ at the heart of our inversion in terms of volumes of the 
It\^o measure $n$ as follows:
$$p(\sqrt{2D\lambda}\,)=(\sqrt{2}/2) \, 
{\la (\varphi_D)(\lambda)
\over n(\{ \hbox{$e\ge 0$ containing a $\lambda$-marked point}\})}\,,\quad \lambda>0\,,
$$
setting  
$\varphi_D(t)={\bf 1}_{\{t>D\}} \,n(\{\hbox{$e\ge 0$ of ${\rm lifetime}(e)>t$}\} )$,
for any real $t>0$.
\vskip8pt
\centerline{\headingfont 5.\quad 
Laplace transforms of the excursion law: principal results}
\vskip.08cm
In this section we identify a principal structure of the excursion law
enforced when working with processes which are restartable twice in the 
sense of Section~2.4. The thrust is to seek giving expression to the 
effects of the two sources of stochasticity featured in Section~3 by 
admitting the use of transform methods in time direction, and to 
identify a framework of sufficient conditions for rendering this 
rigorous; this framework is then to be verified in particular when 
working with Brownian motion in the further development of the paper.  
\ppar
We therefore concentrate now on the Section~2.4 normalized excursion laws 
$h_b^*$ associated with Markov processes $X$ restartable-at-a-stopping-time 
twice (with the processes that result denoted by $X^*$ and $X^{**}=(X^*)^*$ 
respectively), and, for simplicity, assume absolute continuity with respect 
to Lebesgue measure of all random variables under consideration to hold
as well. For the study of the existence of $h_b^*$ we adopt the Laplace 
transform with respect to time, $\la$; see Appendix~A for pertinent 
notation and concepts. 
\goodbreak\medskip
{\bf 5.1\quad Basic Laplace transform:}\quad Our starting point is furnished by 
the representation of $h_b^*$ in terms of the transition densities 
of $X$  as follows:
$$ h_b^*(u,y)=E\big[ {\bf 1}_{\{H_b^*<u\}}
p_X\big(u\!-\!H_b^*, X^*_{H_b^*}, y\big)\big],
$$
for any reals $u>0$ and $y$, where  $p_X(\tau,\xi,dx)$
denotes the density of $X$ conditional on time-$\tau$ information
subject to  $\xi=X_t$. To obtain this representation write in 
Section~2.3 the expectation defining $h_b^*$ as an iterated one by taking 
it conditional on time-$H_b^*$ information first, and then express these 
conditional expectations by time-$H_b^*$ restarting of $X^*$ in terms of 
the transition densities $p_X$. 
\ppar
Two sources of stochasticity are thus identified to determine the structure 
of $h_b^*$, namely $H_b^*$ and $X^*(H_b^*)$. They can be separated on taking 
Laplace transforms with respect to time once their independence is granted, 
and the following result is seen~to~hold~then.
\medskip
{\bf Proposition 5.1:}\quad {\it 
Let the Markov process $X$ be twice restartable and 
assume $H_b^*$ and $X^*(H_b^*)$ to be independent. Then we have 
for any real $y$ the~Laplace~transform~identity:
$$
\la\big(h_b^*(\cdot , y )\big)(z)=E\big[\exp(-zH_b^*)\big]\,
E\big[\la\big(p_X(\cdot, X^*_{H_b^*}, y)\big)(z)\big],$$
for any complex $z$ with $\re(z)\ge z_0$, in the sense of measurable
functions with both sides either $\infty$ or finite.}
\medskip
As a next step we distinguish the effects of the Section~2.3 Cases~I 
and II on $H_b^*$. In (II-1) observe $H_b^*=\delta_0$. In (I) and (II-2)  
a reduction occurs to a time-$0$ and excursion-level-$0$ 
case~I situation by way of the decomposition of $H_b^*$ into 
independent random variables:
$$H_b^*=T_b^*+ H_{0,0}^{**}
$$
which gives expression to the very construction of $H_b^*$ in these cases 
as follows: 
the construction refers to the process $X^{**}$ obtained by restarting of 
$X^*$ at time $T_b^*$, the first passage time of $X^*$ to the level $b$, 
and starting from $T_b^*$ it proceeds by measuring $H_{0,0}^{**}$, the 
achievement time associated with $X^{**}$ by the Section~2.2 construction 
in the case~I situation where $a=0$ and $t=0$ there. Observe that $T_b^*$ 
is equal in law to $T_b$, the first passage time of $X$ to the level $b$, 
and $H_{0,0}^{**}$ is equal in law to $H_0^*$, the Section~2.4 normalized 
achievement time for $X^*$ with $b=0$ there. 
Effects thereof are explained in Sections~5.2~and~5.3. 
\medskip
{\bf 5.2\quad Case I specializations:}\quad First consider the 
Case~I situation as left in Section~5.1. Using the independence 
of the summands in the decomposition of $H_b^*$ derived there, 
the Laplace transform of Proposition~5.1 is here checked 
to become a product of three Laplace transforms as follows.
\medskip
{\bf Proposition 5.2:}\quad {\it Let the Markov process $X$ be twice 
restartable. Assume a Case~I situation 
in which the random variables $H_b^*$ and $X^*(H_b^*)$ are independent 
and the Laplace transforms $\la(h_b^*(\cdot,y))$ for any real $y$ are 
well-defined on the half-plane $\{\re(z)>z_0\}$ within the right-hand 
complex half-plane.
Then, for any complex $z$ with $\re(z)> z_0$,} 
$$  \la\big(h_b^*(\cdot,y)\big)(z)=
E[e^{-zH_0^*}\,]\, E[e^{-zT_b^*}\,]\,
E\big[\la\big(p_X(\cdot, b\!+\!X^*_{H_0^*}, y)\big)(z)\big]\,.$$
\smallskip
{\bf 5.3\quad Case II specializations:}\quad Concentrate on the Case~II 
situation as left in Section~5.1.
Referring to Section~2.3 we hence ask about the excursion of  $X^*$ to 
continue below the level $b>0$ for a period of length at least 
$\delta:=\delta_0$ which is positive and smaller than $D$: $\delta\in(0,D)$. 
The principal structure of the Laplace transform of the
normalized excursion law $h_b^*$ here  mirrors an additional
characteristic feature:  the relative position of $\delta$ to $T_b^*$,
the first  passage time of $X^*$ to the level $b$. This affords a
decomposition of the sample space into the
disjoint subsets $A_1=\{T_b^*\ge \delta\}$ and $A_2=\{T_b^*< \delta\}$,
corresponding to the Subcases II-1 and II-2 of Section~2.4 respectively.
In Proposition~5.1 factors of the Laplace transform thus
decompose into two term sums. A four term decomposition of the Laplace
transform results whose summands are Laplace transforms as
well. We hence obtain four functions $h^*_{b,j,k}$ on
$[0,\infty)\times{\bf R}$, indexed by $j$, $k$ in $\{1,2\}$, by way of
the Laplace transform identities:
$$
\la\big(h_{b,j,k}^*(\cdot , y )\big)(z)
=E\big[{\bf 1}_{A_j}\, \exp(-zH_b^*)\big]\,
E\big[{\bf 1}_{A_k}\,\la\big(p_X(\cdot, X^*_{H_b^*}, y)\big)(z)\big]
,$$
for any complex $z$ with $\re(z)\ge z_0$. The precise decomposition
result is as follows.
\medskip
{\bf Proposition 5.3:}\quad {\it Let the Markov process $X$ be twice 
restartable. Assume a Case~II situation in which the
random variables $H_b^*$ and $X^*(H_b^*)$ are independent and the Laplace
transforms $\la(h_b^*(\cdot,y))$ for any real $y$ are well-defined on the
half-plane $\{\re(z)>z_0\}$ within the right-hand complex half-plane. Then
the functions $h^*_{b,j,k}$ are well-defined and afford the four-term
decomposition:
$$\displaylines{
\hbox{$
h^*_b=\sum\nolimits_{j=1}^2 \sum\nolimits_{k=1}^2 h^*_{b,j,k}
$}\,.\cr 
\noalign{\vskip-3pt Explicitly we here have:}
h^*_{b,1,1}(u,y)=Q(T_b^*>\delta)\,{\bf 1}_{\{u>\delta\}}
E\big[{\bf 1}_{\{ M(X^*)_\delta<b\}}\,p_X\big(u\!-\!\delta,
X^*_\delta, y\big)\big]\,,\cr
h^*_{b,1,2}(u,y)=Q(T_b^*>\delta){\bf 1}_{\{u>\delta\}}
E\big[{\bf 1}_{\{ T_b^*\le\delta\}}\,p_X\big(u\!-\!\delta,
b\!+\!X^*_{H_0^*}, y\big)\big]\,,\cr}$$
for any real $u>0$ and $y$, using  the running maximum
$M(X^*)_\delta=\max\{ X^*_\tau\,|\,\tau\in [0,\delta]\}$. The remaining
two functions are for any real $y$ characterized by the three-factor
Laplace transform identities on $\{\re(z)>z_0\}$ as follows:}
$$\eqalign{
\la\big(h_{b,2,1}^*(\cdot , y )\big)(z)
&=E\big[e^{-zH_0^*}\, \big]\,
E\big[{\bf 1}_{\{T_b^*\le \delta\}}\,e^{-zT_b^*}\,\big]\,
E\big[{\bf 1}_{\{ M(X^*)_\delta<b\}}\,
\la\big(p_X(\cdot, X^*_\delta, y)\big)(z)\big],\cr
\la\big(h_{b,2,2}^*(\cdot , y )\big)(z)
&=E\big[e^{-zH_0^*}\,\big]\,
E\big[{\bf 1}_{\{T_b^*\le \delta\}}\,e^{-zT_b^*}\,\big]\,
E\big[{\bf 1}_{\{ T_b^*\le \delta\}}\,
\la\big(p_X(\cdot, b\!+\!X^*_{H_0^*}, y)\big)(z)\big].\cr}
$$
\smallskip
This result is essentially extracted from [13, Appendix Section 8.3] and
ultimately makes explicit the effects of the Section~5.1 restarting at
$T_b^*$ decomposition of $H_b^*$ once more. Seen in conjunction with
Proposition~5.2, rendering $H_b^*$ explicit thus reduces to
making explicit the law of $H_0^*$ and $X^*(X_0^*)$. We look at this
in the Brownian case.
\medskip
\vskip8pt
\centerline{\headingfont 6.\quad Review: Laplace transforms in the Brownian case}
\vskip.08cm
The general structure of the Laplace transforms of $h_b^*$ has been 
identified in Section~5. We now revert to the Brownian setting of 
Section~3, and as a second step in the proof of the\goodbreak  
Section~3 results 
recall how to make the Section~5.2 and 5.3 Laplace transforms of $h_b^*$
explicit when $X=W^*$ is Brownian motion there. The argument
follows [13, Appendix]; it is  based on the structure theory of the
Brownian meander and the Az\'ema  martingale as developed in [1] and
[2], with another exposition in [14, Section~12.3].
\medskip
{\bf 6.1\quad Brownian motion at normalized achievement time:}\quad
As a first step in obtaining the Laplace
transforms this section reports how the assumptions of the Section~4
results are satisfied in the
Brownian case as follows.
\medskip
{\bf Fact 6.1:}\quad{\it For any real $b$, the random variables $H_b^*$
and $W^*(H_b^*)$ are independent, and for the law of the absolute value of
the latter we have on {\bf R}:}
$$ |W^*_{H_b^*}|(dx)
= {\bf 1}_{(0,\infty)}(x)\, x\exp\big( -\hbox{${1\over 2}$} x^2\big)\, dx
\,.$$
Simplifying notation $B=W^*$ we recall  how this is based
on the {\it Brownian meander\/} which for any real $T>0$ is the process
$m_T$ given by
$$m_T(u)= \big|
W^*\bigl( g_T+u(t\!-\!g_T)\bigr)\big|/{\sqrt{ T\!-\!g_T\, }}, 
\quad u\in[0,1];$$
Here $g_T=\sup\{s\,|\,\hbox{$s\le T$ and $B(s)=0$}\}$ is the last time
before $T$ where $B$ is $0$. The idea is to study  $m_T$
with $T=H_b^*$. Since $\fa_{g_T}\subseteq \fa_{T}$ by
[10, XII (3.2) Lemma, p.~464], it is sufficient to show that 
$m_T$ is independent of $ \fa_{g_T}$ and ${\rm sgn}(B_T)$, and 
has the law indicated independently of $u$. Since $g_T=T g_1$ by 
Brownian scaling, a further reduction occurs to the case $T=1$ where
these independence properties are  taken care of by
the following result.
\medskip
{\bf Fact 6.2:}\quad {\it  With respect to the measure
$Q^*=(|B_1|/E[|B_1|])Q$, the process $m_1$ is a dimension $3$ Bessel process
${\rm BES}(3)$ which is independent of $\fa_{g_1}$.}
\medskip
In [2, \S Th\'eor\`eme, p.~293] or [14, Section~12.3.2]
this is proved as a consequence of a generalized Girsanov
argument; the idea is to thus construct a Brownian motion $\beta$ such 
that $m_1$ satisfies
$m_1(u)=\beta(u)\!+\! \int_0^u (1/m_1(s))\, ds $ as it characterizes
Bessel processes of dimension $3$ started at $0$ at time $0$.
\medskip
As indicated in [2, p.~294] or [14, Step~3, p.~45], the
independence properties of Fact~6.2 entail a relation going
back at least to [7] as follows.
\medskip
{\bf Fact 6.3:}\quad {\it We have
${\rm law}(m_1)(u,x)=\sqrt{\pi/2}\, x^{-1} \, {\rm law}({\rm BES}(3))(1,x)$
for any $u$ in $[0,1]$ and any real $x>0$, where 
${\rm law}({\rm BES}(3))(1,x)=2^{-1/2}\Gamma(3/2)^{-1} x^2 \exp(-x^2/2)$.}
\medskip
{\bf 6.2\quad Case I Laplace transforms:}\quad Keeping the concepts
and the notation of Section~6.1, we  look at the consequences of
the results there for the Laplace transforms of the normalized 
excursion law $h_b^*$ in Case~I, and indicate how a reduction occurs to the 
Laplace transform of the Section~2.4 normalized achievement time $H_0^*$.
Indeed, as a consequence of Fact~6.1 the independence assumptions 
for Proposition~5.2 are satisfied.  The law of Brownian motion at time 
$H_0^*$ as it enters into the Laplace transform representation there is 
explicitly known by Fact~6.1 as well.  Summarizing the computations of 
[11, Sections~8 and 9] following [13, Appendix Sections~8.3.1 and 8.3.2], 
we thus have for any complex $z$ with $\re(z)>0$ the representation:
$$\displaylines{
\la(h_b^*(\cdot, y))(z)= E^*[\exp(-zH_0^*)]\, R_I(y,\sqrt{2Dz}\,),\cr
\noalign{in the sense of measurable functions, where the factor}
\eqalign{
R_I(y,\sqrt{2Dz}\,)&\deff E^*[\exp(-zT_b^*)] E^*[ \la(\chi_{|X_2-y|})(2z)]
\cr &=\exp\big(b^*\sqrt{2Dz})(\sqrt{D}\, G_{\beta(y)})(2Dz)\,,\cr
}\cr}$$
with $b^*=b/\sqrt{D}\le 0$, is well-defined and finite on $\{\re(z)>0\}$.
Here recall the functions $\chi_a$ from Section~A.1, $G_\alpha$ from
Section~A.3, and $\beta$ from Section~3.1.
\goodbreak
\medskip
{\bf 6.3\quad Case II Laplace transforms:}\quad  Keeping the concepts and
the notation of Section~6.1, we  study the Case~II situation. Here a
reasoning analogous to that of Section~6.2 yields the Laplace transform 
representation of  $h_b^*$ of  Proposition~5.3 in the sense of measurable 
functions. A reduction occurs of this law to that of determining two Laplace 
transforms. The computations of 
[11, Section~10] following [13, Appendix Section~8.3.3] yield for the the 
first of these, the Laplace transform of the density
$h_{b,2,1}^*$, the representation:
$$\displaylines{
\la(h_{b,2,1}^*(\cdot, y))(z)= \ndint 0 {\delta}
\exp(-zw)\,E^*[\exp(-zH_0^*)]\, R_{2,1}(y,\sqrt{2Dz}\,)\mu_b(dw)\,,\cr
\noalign{for any complex $z$ with $\re(z)>0$, in the sense of
measurable functions. Here the factor}
\eqalign{
R_{2,1}(y,\sqrt{2Dz}\,)&\deff E^*[{\bf 1}_{\{ T_b^*> \delta\}}\,
\la(\chi_{|X_2-y|})(2z)]
\cr &=\sqrt{D}\ndint {\bf R} {}
\la(\chi_{|x-y|/\sqrt{D}})(2Dz)\varphi_{b,\delta}(x)\, dx\,\cr
}\cr }$$
is well-defined and finite on $\{\re(z)>0\}$, and the remaining concepts are
as follows: from Section~A.1 we have the density $\mu_b$ of $T_b^*$,
the first passage time of Brownian motion $W^*$ to the level $b$, as well as 
the functions $\chi_a$, while the function $\varphi_{b,\delta}$ is
from Section~3.3.
\ppar
For the second Laplace transform in question, the one of the density
$h_{b,2,1}^*$, we have  for any complex $z$ with $\re(z)>0$ the
representation:
$$\displaylines{
\la(h_{b,2,2}^*(\cdot, y))(z)= \ndint 0 {\delta}
\exp(-zw)\,E^*[\exp(-zH_0^*)]\, R_{2,2}(y,\sqrt{2Dz}\,)\mu_b(dw)\,,\cr
\noalign{in the sense of measurable functions, where the factor}
\eqalign{
R_{2,2}(y,\sqrt{2Dz}\,)&\deff E^*[{\bf 1}_{\{ T_b^*\le \delta\}}\,
\la(\chi_{|X_2-y|})(2z)]
\cr &=\sqrt{D}\, Q^*(T_b^*\le \delta)\, G_{\beta(y)}(2Dz)\,\cr}\cr
}$$
is well-defined and finite on $\{\re(z)>0\}$. Here again recall the
Section~A.3 function $G_\alpha$.
\medskip
{\bf 6.4\quad The achievement time Laplace transform:}\enspace 
This section establishes the Laplace transform of the 
normalized achievement time $H_0^*$ by recalling ideas of
[13, Appendix].
The principal finding is that the thus encoded excursion theoretic aspect
of the problem leads to higher transcendental function denominators in
the Laplace transform  as follows:
$$
E[\exp(-zH_0^*)]={1/\Psi(\sqrt{2Dz}\,)}\,, \quad \re(z)>0\,,
$$
with $\Psi$ the Section~A.2 function.
The significance of this result at this stage is that it provides the 
missing factor in the Laplace transforms of Sections~6.2 and 6.3 and 
also shows all Laplace transforms
to be finite~on~$\{\re(z)>0\}$.
\ppar
The equivalent form of this result to be proved is
the following {\it key relation\/}:
$$ E\big[ \exp\big(-\hbox{${1\over 2}$} z^2 H_0^*\big)\big] \,
\Psi(z\sqrt{D\,})=\Psi(0)=1,$$
for any complex $z$ with $\re(z^2)>0$. Taking up the discussion of
Section~6.1, the key relation is based on the {\it Az\'ema
martingale\/} $\mu$. This process is defined in 
terms of the\goodbreak
Brownian meander of Section~6.1 by way of the equality:
$$W^*_u=m_u(1)\, \mu_u\,,\quad u\in [0,\infty);$$
it furnishes a martingale with respect to ${\bf F}^+$, the 
progressive enlargement of the Brownian filtration ${\bf F}$ 
by the sign of $W^*$ whose time-$u$ step is given by
${\bf F}^+(u)=\fa_{g_u}\vee \sigma({\rm sgn}\, W_u^* )$, for any $u\ge 0$.
Taking  stochastic exponentials of this defining equality 
using the independence results of Fact~6.1 at time $1$ there,                                                         
$$
E\big[ \exp\big(zW^*_t\!-\!\hbox{${1\over 2}$}z^2t\big)
\, \big|\, \fh^{\, +}(g_t)\big] =
\exp\big(-\hbox{${1\over 2}$}z^2 t\big)\, \Psi\big(z\, \mu_t\big)\, , $$
for any  real $z>0$. For fixed such $z$, the right-hand side of the 
equality at time $t=H_0^*$ is a martingale by
appealing to forms of the optional stopping theorem.  The expectation
of this martingale stopped at $H_0^*$ is equal to its time-$0$
expectation which is equal to $\Psi(0)$, and hence equals $1$. Use the 
independence of  $H_0$ and $W^*(H_0^*)$ to obtain the  key relation.
\vskip8pt
\centerline{\bf 7.\quad Identification of the Cases I and II inversion problems}
\vskip.1cm
We obtain the results of Section~3 in three steps from those of 
Sections~5 and 6 proceeding by reduction to the situation addressed 
by Theorem~4.1. Summarizing Sections~5 and 6, this section provides 
the first of these steps and identifies the functional relations to 
be considered. Since these transcribe the effects of two sources of 
stochasticity in two principal situations it is two representations 
which thus result as follows.
\medskip
{\bf Lemma 7.1:}\quad {\it For any reals $u>0$ and $y$, we have in Case~I:}
$$
h_b^*(u,y)={1\over 2D}\la^{-1}\Big(
{R_I(y,\sqrt{z}\,) \over \Psi(\sqrt{z}\, )}\Big) \Big({u\over 2D}\Big).$$
\smallskip
{\bf Lemma 7.2:}\quad {\it For any reals $u>0$ and $y$ and with $k$ in
$\{1,2\}$, we have in Case~II:}
$$\eqalign{
h_{b,2,k}^*(u,y)&={1\over 2D}\ndint 0 {\delta}\,  \mu_b(dw)\, \la^{-1}\Big(
\exp\Big(\! -\, {w\over 2D}z\Big)
{R_{2,k}(y,\sqrt{z}\,) \over \Psi(\sqrt{z}\, )}\Big)
\Big({u\!-\!w\over 2D}\Big), \cr
&={1\over 2D}\ndint 0 {\min\{\delta,u\}}\,  \mu_b(dw)\, \la^{-1}\Big(
{R_{2,k}(y,\sqrt{z}\, )\over \Psi(\sqrt{z}\, )}\Big)
\Big({u\!-\!w\over 2D}\Big). \cr }$$
Here Lemma~7.1 summarizes the findings of Section~6.2 and Lemma~7.2 those 
of Section~6.3, after a change of variables in both cases. 
\vskip8pt
\centerline{\headingfont 8.\quad Identification of the auxiliary random variables}
\vskip.08cm
On comparison with Theorem~4.1 the results of Section~7 suggest immediate 
candidates for the functions $R$ governing the auxiliary random variables 
$A_R$ and $B_R$ to be introduced there. As a second step of the inversion 
procedure for the Section~3 results this section therefore identifies these 
functions and verifies their pertinent properties. These functions, it 
might be worth recalling from Section 6, transcribe the effects of two 
sources of stochasticity and their interrelations. In what follows let $y$ 
denote an arbitrary real.
\medskip
{\bf Lemma 8.1:}\quad {\it The function $R_I(y,\sqrt{z}\,)$ on the
right-hand complex half-plane $\{\re(z)>0\}$ is obtained as a Laplace
transform and
$R_I(y,\sqrt{z}\,)=O(z^{-3/2})$ as $z$ tends to $\infty$ there.}
\medskip
{\bf Lemma 8.2:}\quad {\it We have
$\la^{-1}(R_I(y,\sqrt{z}\,))=\sqrt{D}\, \rho_{b^*,\beta(y)}$.}
\medskip
Proving these two results together we have to look at the functions
$R(z)=R_I(y,\sqrt{z}\, )=\sqrt{D}\exp(b^*z)G_{\beta(y)}(z)$ where
$b^*=b/\sqrt{D}\le 0$, and we address their inversion  right away by\goodbreak
replicating the functions  $\rho_{a,c}$ of Section~3.1 with $a=b^*$ and
$c=\beta(y)$. First let  $\beta(y)\le 0$. Setting
$\alpha(y)=-(b^*\!+\!\beta(y))$, then appeal to Remark~A.3 o obtain:
$$
\la^{-1}\Big( {R(z)\over \sqrt{z}}\Big)=
\sqrt{D}\,\la^{-1}\Big( {\exp(-\alpha(y)\sqrt{z}\,)\over \sqrt{z}}\,
\HHfactor_1(z)\Big)
=\sqrt{D}\,\chi_{\alpha(y)}*\hhfactor
$$
using Section~A.1 for inverting the first factor and
Proposition~A.1 for inverting the second one. In the
case $\beta(y)>0$, on the other hand,
$$
\la^{-1}\Big( {R(z)\over \sqrt{z}}\Big)
=\sqrt{D}\,\la^{-1}\Big( {\exp(-|b^*|\sqrt{z}\,)\over \sqrt{z}}\,
G_{\beta(y)}(z)\Big)
=\sqrt{D}\,\chi_{|b^*|}*g_{\beta(y)}$$
now using Proposition~A.2 for inverting the second factor. The proof of
Lemma~8.2 is complete. Fom this discussion the asymptotic behaviour of $R$
near $\infty$ required is immediate if $b^*$ or $\beta(y)$ are not $0$. 
Otherwise, this behaviour is inherited from the one of $N_1$ by way 
of Appendix~B, and the proof of Lemma~8.1 is complete as well. 
\medskip
{\bf Lemma 8.3:}\quad {\it The function $R_{2,1}(y,\sqrt{z}\,)$ on
$\{\re(z)>0\}$ is obtained as a Laplace transform and
$R_{2,1}(y,\sqrt{z}\,)=O(z^{-3/2})$ as $z$ tends to $\infty$ there.}
\medskip
{\bf Lemma 8.4:}\quad {\it We have
$\la^{-1}(R_{2,1}(y,\sqrt{z}\,))=\sqrt{D}\, \rho_{0,\beta(y)}$.}
\medskip
These are special cases of Lemma~8.1 and 8.2 respectively;
from the definitions in Sections~3.3 and 6.3  recall how the
functions $R_{2,1}$ arise from functions $R_I$ with $b^*=0$. 
\medskip
{\bf Lemma 8.5:}\quad {\it The function $R_{2,2}(y,\sqrt{z}\,)$ on
$\{\re(z)>0\}$ is obtained as a Laplace transform and
$R_{2,2}(y,\sqrt{z}\,)=O(z^{-1})$ as $z$ tends to $\infty$ there.}
\medskip
{\bf Lemma 8.6:}\quad {\it We have
$\la^{-1}(R_{2,2}(y,\sqrt{z}\,))(u)=\rho_b (u,y)$ for any real $u>0$.}
\medskip
Setting $R(z)=R_{2,2}(y,\sqrt{z}\,)$ recall
$R(z)=\sqrt{D}\int_{\raise2pt \hbox{$\scriptstyle {\bf R}$}}
\la(\chi_{|x-y|/\sqrt{D}})(z)\varphi_{b,\delta}(x)\, dx $
from Section~6.3 in establishing these two results.
\ppar
{\it Proof of Lemma 8.5\/}\enspace With the Section~A.1  Laplace transforms 
$\la(\chi_a)$ defined on the right-hand half-plane, the function $R$ is defined 
there as well.
To determine its asymptotics there, express $R$ by a direct computation
in terms of the function $\erfc$ as follows:
$$
\sqrt{z}\, R(z)=\sqrt{2\delta D}\,({\sqrt{\pi}/ 2})
\sum\nolimits_{\varepsilon\in\{\pm 1\}} \exp(\delta z \!+\! \varepsilon
y \sqrt{2z}\,)\erfc\big(
(1/\sqrt{2\delta})(\delta\sqrt{2z}\!+\!\varepsilon y)\big)\,.
$$
Using the leading term of the asymptotic expansion of $\erfc$
on the right-hand half-plane (see [9, Section~2.2]) the right-hand
side of this expression is checked to behave like a
scalar multiple of $1/\sqrt{z}$  as  $z$ tends to $\infty$ there.
Hence $R(z)=O(z^{-1})$, as desired. 
\ppar
To exhibit $R$ as a Laplace transform, and to thus complete the proof,  it 
is tempting to argue that its Laplace inverse is obtained by 
inversion under the sign of the defining integral. This turns out to be 
correct proceeding in  two steps as follows. First delete  an 
$\varepsilon$-neighbourhood of $y$ from the domain of integration of 
the defining integral, and effect Laplace inversion by inversion under 
the integral sign. Identify the integrand thus obtained as a 
continuous function bounded by an integrable function independent of $y$, 
and let $\varepsilon$ shrink to $0$. 
\ppar
{\it Proof of Lemma 8.6\/}\enspace 
The last  two step argument in proving Lemma~8.5 extends to the situation 
of Lemma~8.6. Reminding the Section~A.1 function $\varphi_\alpha$, we  here 
obtain
$$
\la^{-1}\Big( {R(z)\over \sqrt{z}}\Big)(\tau)=\sqrt{D}\ndint {\bf R} { }
\varphi_{b,\delta}(x)\, \erfc\Big( {1\over 2}\,
{|x\!-\!y|\over \sqrt{D\tau}}\Big) dx=\rho_b(\tau,y)\,,
$$
where the last equality results noting
$\erfc(z/\sqrt{2D\tau})=2N_{0, D\tau}(z)$.
The proof of Lem\-ma~8.6 is complete as well.
\goodbreak
\vskip8pt
\centerline{\headingfont 9.\quad Proof of the excursion density results}
\vskip.08cm
Summarizing the development up to now, this section establishes the three
convolution sum representation results, Theorems~3.1, 3.2 and 3.3, for the
normalized excursion law and its two distinguished summands respectively.
Referring to Section~7, each of these functions is the inverse of an
explicitly given Laplace transform, and the principal idea is to effect these
inversions analytically by applying Theorem~4.1. Apart from collecting terms
the task at this point therefore reduces to two things. First, to establish
for the numerators $R$ of the Laplace transforms of Section~7 the asymptotic
behaviour required in Theorem~4.1. Second, to identify the Laplace
inverses of these
numerators on division by the complex square root. This we address in
turn in the three sections to follow.
\medskip
{\bf 9.1\quad Proof of Theorem~3.1:} To establish the Theorem~3.1
description of the normalized excursion density $h_b^*$ in Case~I,
we start from the Laplace inversion problem of Lemma~7.1. The idea is
to apply Theorem~4.1  with the functions
$R(z)=R_I(y,\sqrt{z}\, )=\sqrt{D}\exp(b^*z)G_{\beta(y)}(z)$ on
$\{\re(z)>0\}$. With the assumptions for this result satisfied by
Lemma~8.1 indeed,
$$
h^*_b(u,y)={1\over 2D}\sum\nolimits_{n=1}^\infty {(-1)^{n-1}\over
(2\pi)^{n/2}}\, {\bf 1}_{({n\over 2},\infty)}\Big({u\over 2D}\Big)
\la^{-1}\Big( {R(z)\over \sqrt{z}}\Big)* \hhfactor_{n-1}
\Big( {u\over 2D}-{n\over 2}\Big),
$$
where the summation is over the finitely many positive integers
$n$ satisfying $n<u/D$. The Laplace inverses here are taken care of by
Lemma~8.2 to be equal to $\rho_{b^*,\beta(y)}$, and the proof of
Theorem~3.1 is complete.
\medskip
{\bf 9.2\quad Proof of Theorem~3.2:}\quad To establish the Theorem~3.2
description of the summand $h_{b,2,1}^*$  of the excursion law in Case~II,
we start from the inversion problem of Lem\-ma~7.2 with $k=1$ there. Recall
that this asks to integrate a Laplace inverse with respect to  a density, and
therefore we first look at the inversion problem alone. The idea is to
apply for this inversion Theorem~4.1
with the function $R(z)=R_{2,1}(y,\sqrt{z}\,)$. With the assumptions of this
result satisfied by Lemma~8.5 indeed,
$$
\la^{-1}\Big({R(z)\over \Psi(\sqrt{z}\,)}\Big)(\tau)
=\sum\nolimits_{n=1}^\infty {(-1)^{n-1}\over (2\pi)^{n/2}}\,
{\bf 1}_{({n\over 2},\infty)}(\tau)\,
\la^{-1}\Big( {R(z)\over \sqrt{z}}\Big)* 
\hhfactor_{n-1} \Big( \tau\!-{n\over 2}\Big),
$$
for any  $\tau>0$, where only summands to indices $n<2\tau$ are not $0$.
The Laplace inverses here are then taken care of by Lemma~8.6 to be equal
to $\rho_{0,\beta(y)}$, as desired.
\ppar
A technical problem occurs on substitution of these results in
Lemma~7.2. This is because one has to integrate there with respect to 
$\mu_b(dw)$, and the integration variable enters into the
respective number of the summands. Noting that these numbers are smaller
than $u/D$, however, the representation of Theorem~3.2 is checked to follow,
and the proof of this result is complete.
\medskip
{\bf 9.3\quad Proof of Theorem~3.3:}\quad Establishing the Theorem~3.3
description of the summand $h_{b,2,2}^*$ of the excursion law in Case~II
essentially reduces to the problem considered in Section~9.1. In fact,
writing out the second Laplace inversion problem of Lemma~7.2 with
$k=2$ there gives the representation:
$$
h^*_{b,2,2}(u,y)
={1\over 2\sqrt{D}} \ndint 0 {m_u} \la^{-1}\Big({G_{\beta(y)}(z)\over
\Psi(\sqrt{z}\,)}\Big)
\Big( {u\!-\!w\over 2D}\Big)\mu_b(dw)\,,$$
where $m_u=\min\{\delta,u\}$. The Laplace inverses here are those
of Sections~4.1 and 9.1 when $b^*=0$ there, as formalized in
Lemma~8.3 and 8.4.
On inspection, substitution of the latter result's Laplace inverses  yields
the desired representation  on harmonizing the number of summands as in
Section~9.2. The proof of Theorem~3.3 is complete.
\vskip8pt
\centerline{\headingfont 10.\quad Vista}
\vskip.08cm
The paper has introduced a Parisian-style excursion law and determined its
structure in the Brownian case, thus taking up and extending a development 
initiated in  [13]. The characterization of our Brownian excursion law in terms
of sums of independent random variables, however, does not lend itself readily
to actual work with the law as it would be desirable for addressing the 
declared motivation from finance of [13], for example. 
\input xy
\xyoption{curve}
\xyoption{arrow}
%
%
$$
\xy/r.8cm/:
 (0,0)*{\bullet},
 (3.8,0)="T",
 (1,0)*{\scriptstyle \bullet},
(2.0, 0)*{\scriptstyle \bullet},
(3.0,0)*{\scriptstyle |},
(0,-1)*{\scriptstyle -},
(0,-2)*{\scriptstyle -},
(0,-3)*{\scriptstyle -},
(0,-4)*{\scriptstyle -},
 (-.5,0), {\ar (0,0);"T" }, 
 (0,0), {\ar (0,-4.35);(0,0.6)},        
 (1, 0.3)*{\scriptstyle 95},
(2.0, 0.3)*{\scriptstyle 100},
 (3.0,0.3)*{\scriptstyle 105},
 (4.0,0.3)*{\scriptstyle S_t},		
(-.6,.3)*{\scriptstyle\hfill  0},
(-.7,-1)*{\scriptstyle -0.1\,},
(-.7,-2)*{\scriptstyle -0.2},
(-.7,-3)*{\scriptstyle -0.3},
(-.7,-4)*{\scriptstyle -0.4},
(0,0),{(1.0,-4.002543342);
(3.0, -2.544991699839)
\curve{~**\dir{-}
&(1.0,-4.002543342833239)   
&(2.0,-3.216761851275)    
&(3.0,-2.544991699839)    
}},
(0,0),{(1.1, -0.84845716048148);
(3.0, -0.407191452132)
\curve{~**\dir{o}
&(1.1, -0.84845716048148)   
&(2.0, -0.6156599307908)      
&(3.0, -0.407191452132)  
}},
\endxy
$$
\vskip-4pt
\centerline{{\eightbf Figure 10.1.}\enspace \ninerm 
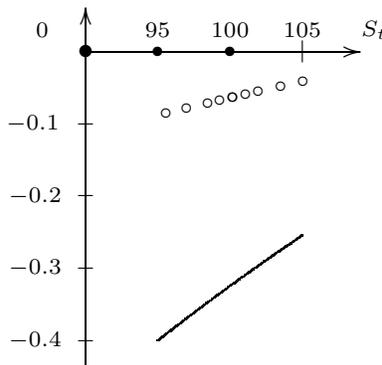
Comparison of Deltas of Down-and-In call [-] and  
Parisian Down-and-In call [o]}
\medskip
Two methods for 
explicit handling of laws of sums of independent random variables are 
thus developed in the companion [12] to the present paper. Anticipating
here some of the results there, for the  Parisian barrier options 
proposed in [13] the methods are found to furnish effective and 
stable ways not just for  valuation of the  but also for hedging; 
this in particular so in 
situations where standard barrier options build up large Deltas `near
their barriers'; a typical example is featured in Figure~10.1 where 
Parisian barrier options permit a reduction of Deltas by some  factor $5$.
We have thus come full circle, with the results of the present paper 
instrumental for this.
\vskip8pt
\centerline{\headingfont Appendix A.\quad
Laplace transform pairs}
\vskip.08cm
This appendix collects pertinent Laplace transform  pairs.
Here the Laplace transform is the linear operator $\la$ on
the continuous functions of exponential type on $(0,\infty)$,
the positive reals, defined as follows: it associates with
any such function $f$ the function
$\la(f)$~given~by:
$$\la(f)(z)=\ndint 0 {\infty} \exp(-zu)f(u)\, du\,,$$
for any complex $z$ in a half-plane contained sufficiently deep within the
right-hand complex half-plane $\{ z|\re(z)>0\}$. The maps $\la(f)$  are
analytic on such half-planes, and the operator $\la$ is an injection
with inverse $\la^{-1}$, the inverse Laplace transform; see
[3] or [4] for more detail. We moreover work with the principal branch
of the complex logarithm on ${\bf C}\setminus(-\infty,0]$, the complex
plane {\bf C} cut along the non-positive reals $(-\infty,0]$.
\medskip
{\bf A.1}\quad   The following
three standard Laplace transforms on
$\{z|\re(z)>0\}$ from [4, Beispiel 8, p.~50f]
originate with the heat equation:
$$\leqalignno{
\la \bigl( \psi_\alpha\bigr)(z)
=\exp(-\alpha \sqrt{z }\,)
\quad &\hbox{where}\quad \psi_\alpha(u)=
{\alpha\over 2\sqrt{\pi u^3 }}\, \exp\Big(-{\alpha^2\over 4u}\Big),
&\cr
\la  \bigl( \chi_\alpha\bigr)(z)={\exp(-\alpha \sqrt{z }\,) \over
\sqrt{z\ }}
\quad &\hbox{where}\quad \chi_\alpha(u)=
{1\over \sqrt{\pi u}}\,
\exp\Big(-{\alpha^2\over 4u}\Big),
\cr
\la  \bigl( \varphi_\alpha\bigr)(z)={\exp(-\alpha \sqrt{z }\,) \over z}
\quad &\hbox{where}\quad \varphi_\alpha(u)=
\erfc\Big({\alpha\over 2\sqrt{u}}\Big),
&\cr
}$$
for any real $u>0$. Here $\erfc(\xi)=(2/\sqrt{\pi})
\int_{\raise 3pt\hbox{$\scriptstyle [\xi,\infty)$}} \exp(-x^2)\, dx$,
for any complex  $\xi$, is the complementary error function, and $\alpha$
is any complex with $|\arg(\alpha)|\le {\pi/4}$ such that
$\re(\alpha)>0$ for~$\psi_\alpha$. Hence the law $\mu_b$ of
the first passage time of Brownian motion to the  level $b$ is given
on the real line {\bf R} by:
$ \mu_b(dw)=\psi_{b\sqrt{2}}(w)\, dw$.   
\medskip
{\bf A.2}\quad As the
first of two sets of functions to be considered define the functions
$\HHfactor_n$ on ${\bf C}\setminus(-\infty,0]$ for any integer $n\ge 0$  by:
$$ \HHfactor_n=\HHfactor_1^n\quad\hbox{where}\quad 
\HHfactor_1(z)=\Psi(-\sqrt{z}\, )/\sqrt{z}\,,
$$
for any $z$ in  ${\bf C}\setminus(-\infty,0]$. Following [16],
the function $\Psi$ here is the generalization of
the normal distribution given by the integral:
$$\Psi(w)=\ndint 0 {\infty}
x\exp\big(-\hbox{${1\over 2}$}x^2\!+\!wx)\, dx\,,\quad w\in{\bf C}\,.
$$
The second set is furnished by the functions $\hhfactor_n=\hhfactor^{*(n)}$ 
on $[0,\infty)$ already used in Section~3. Let  $\hhfactor_0=\hhfactor^{*(0)}$ 
be the Dirac delta function at $0$ and define $\hhfactor_n$ for any integer 
$n\ge 1 $ as an $n$-fold convolution on~$[0,\infty)$~by:
$$
\hhfactor_n=\hhfactor^{*(n)} \quad \hbox{where}\quad
\hhfactor(u)=
(2/\sqrt{\pi }\, )\, \sqrt{u\, }/ (2u\!+\!1)
\, ,\quad
u\in[0,\infty). $$
{\bf Proposition A.1:}\quad{\it We have $\la(\hhfactor_n)=\HHfactor_n$ 
on $\{\re(z)>0\}$, for any integer $n\ge 0$.}
\medskip
\medskip
{\bf A.3}\quad This section
concentrates on generalizations of the Section~A.2 function $\Psi$. For
any real $\alpha$
these are the functions $F_\alpha$ and $G_\alpha$ on
${\bf C}\setminus(-\infty,0]$ given by:
$$\displaylines{
F_\alpha(z)=\ndint 0 {\infty}
x\exp\big(-\hbox{${1\over 2}$} x^2\!-|\alpha\!-\!x|\sqrt{z}\,\big)\,
dx\,, \cr
\noalign{\vskip-2pt and\vskip-3pt}
G_\alpha(z)={F_\alpha(z)/\sqrt{z}}\,.\cr }
$$
In terms of the Section~A.2 function $h$ and the Section~A.1
complementary error function $\erfc$
moreover define the function $g_\alpha$ on $[0,\infty)$ by:
$$
g_\alpha(u)=\exp\Big( \!-{\alpha^2/2\over 1\!+\!2u}\Big)
\Big\{ h(u)\exp\Big( \!-{\alpha^2\over 4u(1\!+\!2u)}\Big) +
{\alpha\over 1\!+\!2u}
\erfc\Big( {\alpha \over \sqrt{4u(1\!+\!2u)} }\Big)\Big\},
$$
for any $u\ge 0$. Then, in particular all $G_\alpha$ are
Laplace transforms as follows.
\medskip
{\bf Proposition A.2:}\quad {\it We have $\la(g_\alpha)=G_\alpha$ on
$\{\re(z)>0\}$.}
\medskip
{\bf Remark A.3:}\quad If $\alpha\le 0$ we have
$F_\alpha(z)=\exp(\alpha \sqrt{z}\,)\Psi(-\sqrt{z}\,)$.
\vskip8pt
\centerline{\headingfont Appendix B.\quad Further
properties of the function $ \Psi$}
\vskip.1cm
This appendix  develops pertinent
properties of the  function $\Psi$ which from Section~A.2 for
any complex number $w$ is given by the integral:
$ \Psi(w)=\ndint 0 \infty  x\exp\big( -\hbox{$1\over 2$}x^2
\!+\!wx\big)\, dx$.
\ppar
Developing the linear exponential factor of the integrand of $\Psi$ in
its series and integrating the resulting series term by term, yields the
following      {\it  series expansion}:
$$ \Psi(w)=\sum\nolimits_{n=0}^\infty a_n \, w^n \quad\hbox{where}\quad
a_n=({2^{n/2}/ n!})\, 
\Gamma\big(  \hbox{$1\over 2$}(n\!+\!2)\big) . 
$$
This series is absolutely convergent for any complex number $w$, and its
convergence is uniform on compact sets.
As a first appplication
it yields the {\it key identity\/}:
$$\Psi(w)=\Psi(-w)+\sqrt{2\pi \, }\, w\,
\exp\big(\hbox{$1\over 2$} w^2\big)$$
which connects the values of $\Psi$ on the right-hand half-plane with those
on the left-hand half-plane and has been noted in an equivalent form  
in establishing [1, Proposition 1 point~2), p.~94]; the identity can also 
be obtained by partial integration of the defining integrals for $\Psi(w)$ 
and $\Psi(-w)$ respectively by way of the partial integration identity:
$$\Psi(-\sqrt{2}w)=1\!-\! \sqrt{ \pi}\, w\exp(w^2)\erfc(w).$$
This identity is the basis of the {\it leading term expansion\/}
for $\Psi$ on the left hand half plane:
$$ \Psi(-w)={1/ w^2 }+\!R_2(w)\quad\hbox{where}\quad
|R_2(w)|\le {6/ |w|^4}\,,
$$
for any complex $w$ with $\re(w)>0$, which
is a special case of a general uniform
asymptotic expansion of $\Psi$ on the left-hand half-plane.
\goodbreak\bigskip
{\paragraphfont References:}
\medskip
{\ninerm
\baselineskip9.75pt
\parein{20pt}{[1]} J. Az\'ema, M. Yor: \'Etude d'une martingale remarquable,
{\nineit S\'em. Proba. XXIII, LNMS} {\ninebf 1372}, 88--130, Heidelberg:
Springer 1989.
\paraus\vskip1.2pt
\parein{20pt}{[2]} J. Az\'ema, M. Yor: Sur les z\'eros des martingales
continues,
{\nineit S\'em. Proba. XXVI, LNMS} {\ninebf 1526}, 248--306, Heidelberg:
Springer 1992.
\paraus\vskip1.2pt
\parein{20pt}{[3]} R. Beals: {\nineit Advanced mathematical
analysis\/}, New York: Springer, 1973.
\paraus\vskip1.2pt
\parein{20pt}{[4]}
G. Doetsch: {\nineit Handbuch deso that we have the Laplace transform
representation there.r Laplace
Transformation\/} I, Basel: Birkh\"auser, 1971.
\paraus\vskip1.2pt
\parein{20pt}{[5]} 
R. Loeffen, I. Czarna, Z. Palmowski: Parisian ruin probabilities for 
spectrally negative L\'evy processes, arXiv: 1102.4055v1 (2011).
\paraus\vskip1.2pt
\parein{20pt}{[6]}
L. Gauthier: Excursion height- and length-related stopping
times, and applications to finance, {\nineit Adv. Appl. Prob.\/}
{\ninebf 34} (2002), 846--868.
\paraus\vskip1.2pt
\parein{20pt}{[7]} J.P. Imhof: Density factorization for Brownian motion
and the three dimensional Bessel processes and applications, {\nineit
J. Appl. Prob.\/} {\ninebf 21} (1984), 500--510.
\paraus\vskip1.2pt
\parein{20pt}{[8]} I. Karatzas, S. Shreve: {\nineit Brownian motion
and stochastic calculus\/} 2nd ed., New York: Springer, 1991.
\paraus\vskip1.2pt
\parein{20pt}{[9]}
N.N. Lebedev:  {\nineit Special functions and
their applications\/}, New York: Dover, 1972.
\paraus\vskip1.2pt
\parein{20pt}{[10]}
D. Revuz, M.Yor: {\nineit Continuous martingales and Brownian motion\/}
2nd ed., Heidelberg: Sprin\-ger,  1994.
\paraus\vskip1.2pt
\parein{20pt}{[11]}
M. Schr\"oder: Brownian excursions and Parisian barrier options:
a note, {\nineit J. Appl. Prob.\/} {\ninebf 40} (2003), 855--864.
\paraus\vskip1.2pt
\parein{20pt}{[12]}
M. Schr\"oder: On a Brownian excursion law, II: Three methods for 
handling convolutions of probability laws.
\paraus\vskip1.2pt
\parein{20pt}{[13]}
M. Yor, M. Jeanblanc--Picqu\'e, M. Chesney : Brownian excursions
and Parisian barrier options, {\nineit Adv. Appl. Prob.\/}
{\ninebf 29} (1997), 165--184.
\paraus\vskip1.2pt
\parein{20pt}{[14]} M. Yor: {\nineit Some aspects of Brownian motion,
Part II\/}, Basel et al.: Birkh\"auser, 1997.
\paraus\vskip1.2pt
}
\medskip
\centerline{\ninecaps Author's address: Keplerstrasse 30, D-69469
Weinheim (Bergstrasse), Germany}
\bye